\documentclass[11pt]{article}

\usepackage[letterpaper,margin=1in]{geometry}
\usepackage[T1]{fontenc}
\usepackage{amsmath,amssymb,amsthm,mathtools,bm,mathrsfs,graphicx,url,microtype,setspace,needspace,xcolor}
\usepackage{eso-pic}
\usepackage{titling}
\usepackage[hidelinks]{hyperref}
\newcommand{\DriftNominalReturn}{99.888}
\newcommand{\DriftCommonReturn}{99.776}
\newcommand{\DriftRobustReturn}{99.761}
\newcommand{\DriftNominalCommonRegret}{1.386}
\newcommand{\DriftCommonCommonRegret}{0.438}
\newcommand{\DriftRobustCommonRegret}{0.823}
\newcommand{\DriftNominalRegret}{2.153}
\newcommand{\DriftCommonRegret}{2.258}
\newcommand{\DriftRobustRegret}{1.019}

\predate{\begin{center}\small}\postdate{\par\end{center}}
\newtheoremstyle{upright}{0.45\baselineskip}{0.45\baselineskip}
  {\normalfont}{}{\itshape}{:}{0.5em}{}
\theoremstyle{upright}
\newtheorem{theorem}{Theorem}

\newtheorem{remark}{Remark}
\AddToHook{env/theorem/before}{\Needspace{6\baselineskip}}
\AddToHook{env/thebibliography/begin}{%
  \small\setstretch{1.02}\setlength{\itemsep}{0.2em}\setlength{\parskip}{0pt}%
}


\newenvironment{drproof}{\par\noindent\hspace*{1em}\textit{Proof: }}
 {\hfill$\blacksquare$\par}
\newcommand{\R}{\mathbb R}
\newcommand{\E}{\mathbb E}
\newcommand{\tr}{\operatorname{tr}}
\newcommand{\col}{\operatorname{col}}
\newcommand{\diag}{\operatorname{diag}}
\newcommand{\cH}{\mathcal H}
\newcommand{\cW}{\mathcal W}

\title{\bfseries
Regret-Optimal Linear-Quadratic Control under Distributional Drift
}

\author{
Lukas-Benedikt Fiechtner\textsuperscript{1}\qquad Jose Blanchet\textsuperscript{1,2}\\[0.6em]
{\small\textsuperscript{1}Institute for Computational and Mathematical Engineering}\\
{\small\textsuperscript{2}Department of Management Science and Engineering}\\
{\small Stanford University, Stanford, CA 94305, USA}
}
\date{}

\hypersetup{
 pdftitle={Regret-Optimal Linear-Quadratic Control under Distributional Drift},
 pdfauthor={Lukas-Benedikt Fiechtner and Jose Blanchet}
}

\begin{document}

\maketitle
\AddToShipoutPictureFG*{%
  \AtPageLowerLeft{%
    \hspace*{\dimexpr1in+\oddsidemargin\relax}%
    \raisebox{18pt}[0pt][0pt]{%
      \parbox[b]{\textwidth}{\normalfont\footnotesize
        Submitted to the 2027 American Control Conference (ACC); under review.\par
        This work has been submitted to the IEEE for possible publication.
        Copyright may be transferred without notice, after which this version may no longer be accessible.}%
    }%
  }%
}
\thispagestyle{empty}
\pagestyle{plain}

\begin{abstract}
We study finite-horizon LQR with independent disturbances whose laws
may drift over time. We address uncertainty in this evolution through
ex-ante distributionally robust regret optimization (DRRO), which
minimizes worst-case excess expected cost relative to the optimal causal
controller that knows the law sequence. To specify the admissible law
sequences, we introduce a temporally coupled moment ambiguity set with a Gelbrich kinetic-energy
budget and a second budget limiting departures from reference
moments and persistent cumulative mismatch.
Despite the coupling of moments across time, affine
disturbance-feedback synthesis reduces exactly to a semidefinite
program. An optimal policy uses past disturbances to predict the
correction to nominal LQR determined by current and future disturbance
means. The SDP solution also yields worst-case moment paths and
corresponding laws. For standard distributionally robust optimization
(DRO), the analogous SDP relaxation can have a duality gap.
Portfolio liquidation experiments show substantially
lower worst-case regret over the drifting ambiguity set than common-law
DRRO and nominal control.
\end{abstract}

\section{Introduction}\label{sec:drift_introduction}

A standard stochastic LQR model assumes i.i.d. disturbances with a
fixed law. However, the actual disturbance distribution encountered during deployment could be
different due to estimation error or distribution shift.
Distributionally robust optimization (DRO) addresses this uncertainty by placing the unknown
law in an ambiguity set of plausible distributions. The controller
then chooses a policy to minimize expected cost in a minimax game
against a fictitious adversary that selects the worst law in the set.

In a standard dynamic DRO formulation, nature may choose a transition
kernel anew at each stage. State--action rectangularity permits independent
choices across states and actions, enabling robust dynamic programming
\cite{iyengar2005robust,nilim2005robust}, but can be pessimistic because
it allows nature to adapt its choice of law to time, state, and action,
which may be an unrealistic modeling assumption. State rectangularity permits
coupling across actions at each state while retaining independent choices
across states and stages \cite{wiesemann_kuhn_rustem_2013}.
In Wasserstein DRO control, state-dependent disturbance laws chosen at each
stage yield dynamic programming \cite{yang_wasserstein_2021} and, with a soft
Wasserstein penalty in a linear-quadratic model, a Riccati recursion
\cite{kim_yang_2023}.
Commitment to a stationary transition kernel can reduce conservatism:
the transition law may depend on the state and action but must remain
unchanged over time for each pair. Nonrectangular ambiguity permits joint
constraints across pairs. For finite-state, finite-action Markov decision
processes (MDPs), Li et al. \cite{li_kuhn_sutter_2026} develop policy-gradient
methods over stationary policies under nonrectangular ambiguity and
infinite-horizon discounted cost. Wang and Si \cite{wang_si_2026} show that
history-dependent policies achieving sublinear expected regret relative to
every admissible kernel's optimal long-run gain are robust-optimal for
average reward.
In linear-quadratic control, nature can instead choose the joint law of
exogenous noises in advance, without state- or action-dependent selection.
Ta\c{s}kesen et al. \cite{NEURIPS2023_3b7a66b2} prove optimality of linear
policies in a partially observed LQG model with independent, zero-mean noise
laws that may vary over time within stagewise Wasserstein balls around
Gaussian reference laws. Nyl\"of et al. \cite{nylof2026correlated} establish
affine optimality for temporally correlated noise under a Wasserstein ball
around an elliptical joint reference law.

However, merely restricting nature's choices does not remove the conservatism inherent
in optimizing worst-case cost. A DRO policy can sacrifice substantial
performance under favorable laws for a small improvement under the
worst law. Ex-ante distributionally robust regret optimization (DRRO)
addresses this issue by measuring, for each candidate law, the
difference between a policy's expected cost and the minimum expected
cost achievable by a causal controller that knows that law. It then
minimizes the worst such regret over the ambiguity set. Work on
ex-ante DRRO has largely concerned static problems with a single
decision \cite{perakis_regret_2008,agarwal_minimax_2022,cho_wasserstein_2024,fiechtner2025wasserstein}.
To our knowledge, our earlier paper \cite{fiechtner2026common} gave the first exact semidefinite
programming (SDP) formulation of multistage ex-ante DRRO for stochastic
control, despite the general hardness of static DRRO
\cite{fiechtner2025wasserstein}. The formulation optimizes affine
disturbance-feedback policies for LQR under \emph{common stage-law
ambiguity}: disturbances are independent draws from one unknown law,
with admissible means and covariances specified by a Gelbrich ball
centered at reference moments. Past disturbances then inform future decisions, and an
optimal affine controller corrects the nominal policy using their
sample mean.

Here we extend this model to environments in which the disturbance law
can change gradually over time. We introduce \emph{drifting stage-law
ambiguity}: nature chooses independent stage laws whose moment path
is coupled across time. A kinetic-energy budget bounds the sum of
weighted squared Gelbrich distances between consecutive moment pairs,
starting from the reference moments. It provides a moment-based
surrogate for the energy of a curve in Wasserstein space, exact for
elliptical laws with a common generator. A second budget limits the
sum of weighted squared distances from the reference moments and penalizes cumulative
moment mismatch, discouraging sustained one-sided departures.
We derive an exact SDP for
ex-ante DRRO over affine disturbance-feedback policies and recover
worst-case laws from its solution. Surprisingly, the analogous two-budget
DRO relaxation need not be exact. An optimal affine policy adjusts
nominal LQR using past disturbance deviations to predict the weighted
combination of current and future means entering the known-law optimal
action. Its weights may differ across observation times,
balancing past information against possible changes in the law, with
prediction errors weighted by their effect on control performance.
In portfolio liquidation experiments, the proposed \emph{Drifting
DRRO} controller achieves substantially lower worst-case regret over
the drifting ambiguity set than both \emph{Common DRRO} and the
nominal controller.

An alternative regret-based approach is ex-post DRRO, which minimizes
the supremum, over admissible disturbance laws, of expected pathwise
regret against a clairvoyant controller with access to the entire
realized disturbance trajectory. Exact finite-horizon convex
reformulations over linear or affine disturbance-feedback policies
are available under Wasserstein and moment ambiguity
\cite{taha_distributionally_2023,taha_bitar_2026}.
Our comparator instead knows the sequence of laws and remains causal
with respect to disturbance realizations, isolating the cost of not
knowing the disturbance model.

Finite-horizon minimax regret has also been studied with a causal
comparator that knows an unknown constant drift
\cite{pmlr-v168-gurevich22a} or scalar dynamics coefficient
\cite{Carruth2023OptimalAC}. Nonstochastic online control studies
bounded adversarial disturbance sequences and gives regret bounds
relative to a fixed class of stabilizing controllers
\cite{agarwal_online_2019}. Foster and Simchowitz
\cite{foster_simchowitz_2020} obtain polylogarithmic regret for known linear
dynamics and known, fixed quadratic costs in this setting. Our objective is exact finite-horizon
minimax synthesis under distributional ambiguity. Finally, Anderson
and Keehan \cite{anderson_keehan_2025} use penalties on consecutive
Wasserstein distances to estimate evolving distributions, without a
control objective.

\emph{Notation.} $\mathbb S^k$ denotes real symmetric $k\times k$ matrices,
with $\succeq$ and $\succ$ denoting the positive semidefinite and definite orders.
We write $I_k$ for the identity, $\tr$ for trace, and
$\lambda_{\min},\lambda_{\max}$ for extreme eigenvalues.
We use the Euclidean norm $\|\cdot\|$, the weighted norm
$\|v\|_M^2:=v^\top Mv$ for $M\succ0$, and the Frobenius inner product
$\langle A,B\rangle:=\tr(A^\top B)$, with $\|A\|_F^2:=\langle A,A\rangle$.
The operators $\col,\diag$ denote vertical stacking and block-diagonal
assembly. $\mathcal P_2(\R^k)$ denotes laws on $\R^k$ with finite second moments.
Independent product laws are written $\bigotimes_t\mathbb P_t$, or
$\mathbb P^{\otimes T}$ for $T$ identical marginals.

\section{Problem Formulation}\label{sec:drift_problem}

\subsection{Dynamics and Regret Objective}\label{subsec:drift_objective}

We consider a linear-quadratic regulator with dynamics
\begin{equation}\label{eq:drift_dynamics}
 x_{t+1}=A_tx_t+B_tu_t+\Xi_tw_t,\quad t=0,\ldots,T-1,
\end{equation}
where $x_t\in\R^n$ is the state, $u_t\in\R^m$ the control, and $w_t\in\R^d$
the disturbance. The matrices are of appropriate dimensions. The initial state is deterministic and known to the controller.
Performance is measured by the quadratic cost
\begin{equation}\label{eq:drift_cost}
 \begin{aligned}
 L={}&\sum_{t=0}^{T-1}
 \begin{pmatrix}x_t\\u_t\end{pmatrix}^{\!\top}
 \begin{pmatrix}Q_t&D_t\\D_t^\top&R_t\end{pmatrix}
 \begin{pmatrix}x_t\\u_t\end{pmatrix}+x_T^\top Q_Tx_T
 \end{aligned}
\end{equation}
with $Q_t$ and $R_t$ symmetric. Indefinite costs are
allowed, subject to the positive-definiteness assumption on the
Riccati matrices $M_t$ defined later in \eqref{eq:drift_M}.

For a joint law $\mathbb P$ of the disturbance sequence
$\mathbf w=(w_0,\ldots,w_{T-1})$, define the expected total cost of a
policy $\pi$ and the optimal expected cost over causal policies by
\begin{equation}\label{eq:drift_values}
 J(\pi;\mathbb P):=\E_{\mathbb P}^{\pi}[L],\qquad
 J^\star(\mathbb P):=\inf_{\tilde\pi\in\Pi_{\mathrm c}}
              J(\tilde\pi;\mathbb P)
\end{equation}
where $\Pi_{\mathrm c}$ is the causal policy class defined later in
Section~\ref{subsec:drift_policies}. The \emph{ex-ante regret} and the
worst-case regret over an uncertainty set $\mathcal P$ of joint laws
of $\mathbf w$ are defined as
\begin{equation}\label{eq:drift_worst_case_regret}
 \mathcal R(\pi;\mathbb P):=J(\pi;\mathbb P)-J^\star(\mathbb P),\quad
 \mathfrak R(\pi):=\sup_{\mathbb P\in\mathcal P}
                         \mathcal R(\pi;\mathbb P).
\end{equation}
Using the affine policy class $\Pi_{\mathrm{aff}}^w$ defined later in
Section~\ref{subsec:drift_policies}, our distributionally robust regret
objective is
\begin{equation*}\tag{DRRO}\label{eq:drift_drro}
 \inf_{\pi\in\Pi_{\mathrm{aff}}^w}\mathfrak R(\pi).
\end{equation*}

\subsection{Temporally Coupled Stage-Law Ambiguity}\label{subsec:drift_ambiguity}

We allow independent disturbances whose stage laws
may drift over time starting from a \emph{reference} law $\hat{\mathbb{P}}$, whose mean
and covariance matrix we denote by $\hat{\mu}$ and $\hat{\Sigma}\succ0$.
More specifically, we take $\mathcal P=\mathcal P_{\mathrm{drift}}$, where
\begin{equation}\label{eq:drift_law_set}
 \mathcal P_{\mathrm{drift}}:=
 \left\{\bigotimes_{t=0}^{T-1}\mathbb P_t\;\middle|\;
 \mathbb P_t\in\mathcal P_2(\R^d),\,
 (\boldsymbol\mu,\boldsymbol\Sigma)\in\cW_{\mathrm{drift}}\right\},
\end{equation}
where $\boldsymbol\mu=\col(\mu_0,\ldots,\mu_{T-1})$ and
$\boldsymbol\Sigma=(\Sigma_0,\ldots,\Sigma_{T-1})$ collect the means
and covariances of the stage laws $\mathbb P_t$.
We also set $\mu_{-1}=\hat\mu$ and $\Sigma_{-1}=\hat\Sigma$.
The moment set $\cW_{\mathrm{drift}}$ is given by
\begin{equation}\label{eq:drift_moment_set}
 \cW_{\mathrm{drift}}:=
 \left\{(\boldsymbol\mu,\boldsymbol\Sigma)\;\middle|\;
 \Sigma_t\succeq0\quad(0\le t<T),\quad
 \mathcal E(\boldsymbol\mu,\boldsymbol\Sigma)\le\delta_{\mathcal E}^2,\,
 \mathcal L(\boldsymbol\mu,\boldsymbol\Sigma)\le\delta_{\mathcal L}^2
 \right\}.
\end{equation}
To specify $\mathcal E$ and $\mathcal L$, we recall the
\emph{Gelbrich distance}
\begin{equation}\label{eq:drift_gelbrich}
 \mathsf G^2\bigl((\mu,S),(\nu,U)\bigr)
 =\|\mu-\nu\|^2+\mathsf B^2(S,U),
\end{equation}
with the \emph{Bures distance} $\mathsf B$ defined by
\begin{equation}\label{eq:drift_bures}
 \mathsf B^2(S,U)=\tr\!\left(S+U
              -2(U^{1/2}SU^{1/2})^{1/2}\right).
\end{equation}
Now choose weights $a_t,b_t>0$ for $0\le t<T$ and $\eta\ge0$,
and budget radii $\delta_{\mathcal E},\delta_{\mathcal L}>0$.
We define the \emph{kinetic-energy} functional $\mathcal E(\boldsymbol\mu,\boldsymbol\Sigma)$ as
\begin{equation}\label{eq:drift_energy}
 \mathcal E(\boldsymbol\mu,\boldsymbol\Sigma)
 :=\sum_{t=0}^{T-1}a_t\,
 \mathsf G^2\bigl((\mu_t,\Sigma_t),(\mu_{t-1},\Sigma_{t-1})\bigr)
\end{equation}
and the \emph{reach-and-persistence} functional $\mathcal L(\boldsymbol\mu,\boldsymbol\Sigma)$ as
\begin{subequations}\label{eq:drift_reach}
\begin{equation}
 \mathcal L(\boldsymbol\mu,\boldsymbol\Sigma)
 :=\mathcal L_1(\boldsymbol\mu,\boldsymbol\Sigma)
   +\eta\mathcal L_2(\boldsymbol\mu,\boldsymbol\Sigma),
\end{equation}
with the \emph{reach} term
\begin{equation}\label{eq:drift_reach_component}
 \mathcal L_1(\boldsymbol\mu,\boldsymbol\Sigma)
 :=\sum_{t=0}^{T-1}b_t\,
       \mathsf G^2\bigl((\mu_t,\Sigma_t),(\hat\mu,\hat\Sigma)\bigr)
\end{equation}
and the \emph{persistence} term
\begin{equation}\label{eq:drift_persistence_component}
 \mathcal L_2(\boldsymbol\mu,\boldsymbol\Sigma):=
 \sum_{t=0}^{T-1}\mathsf G^2\!\Bigl(
       \bigl(\sum_{s=0}^{t}\mu_s,\sum_{s=0}^{t}\Sigma_s\bigr),
       \bigl((t+1)\hat\mu,(t+1)\hat\Sigma\bigr)\Bigr).
\end{equation}
\end{subequations}

The constraints have the following interpretation. A natural
modeling assumption is that the true disturbance laws drift from the
reference law along a curve in Wasserstein space with bounded kinetic
energy, limiting how rapidly the distributions can change.
In the $2$-Wasserstein metric $W_2$, an absolutely continuous
curve $s\mapsto\mathbb Q_s\in\mathcal P_2(\R^d)$ on $[0,\tau]$, starting at
$\mathbb Q_0=\hat{\mathbb P}$, has kinetic energy
\begin{equation}\label{eq:drift_kinetic_energy}
 \operatorname{KE}(\mathbb Q):=\int_0^\tau|\dot{\mathbb Q}_s|_{W_2}^2\,ds,\;
 |\dot{\mathbb Q}_s|_{W_2}:=\lim_{h\to0}
       \frac{W_2(\mathbb Q_{s+h},\mathbb Q_s)}{|h|},
\end{equation}
where the metric derivative $|\dot{\mathbb Q}_s|_{W_2}$ exists almost everywhere
\cite{ambrosio2008gradient}. Interpolating consecutive laws
$\mathbb P_{-1}:=\hat{\mathbb P},\mathbb P_0,\ldots,\mathbb P_{T-1}$
by constant-speed Wasserstein geodesics over intervals of length
$\Delta_t>0$ gives energy $\sum_{t=0}^{T-1}\frac{W_2^2(\mathbb P_t,\mathbb P_{t-1})}{\Delta_t}$.
Gelbrich distance lower-bounds the Wasserstein distance between laws
with those moments, with equality for elliptical distributions sharing
a common generator, including Gaussian pairs
\cite{gelbrich1990}. With $a_t=1/\Delta_t$, the kinetic-energy constraint in
\eqref{eq:drift_moment_set} bounds a moment-based surrogate of the energy
in \eqref{eq:drift_kinetic_energy}.
A kinetic-energy cap still permits continued drift away from the
reference law. The reach term $\mathcal L_1$ limits such drift by bounding each
stage's moment distance from the reference. The laws may also remain
displaced after the drift stops: a common stage law with nonnominal
moments contributes to $\mathcal E$ only at the first time step,
since its law remains constant afterwards. Such persistent departures
may be unrealistic in some applications. The reach term $\mathcal L_1$
penalizes persistence at every displaced
stage. Additionally, $\mathcal L_2$ penalizes
\emph{cumulative moment mismatch}. At each stage $t$, it compares the
mean and covariance of $\sum_{s=0}^{t}w_s$ with those of $t+1$
independent draws from the reference law.
The cumulative penalty therefore favors compensating shifts and can
continue to accrue after a return to the reference if earlier deviations
remain unoffset.
\subsection{Policies}\label{subsec:drift_policies}

At time $t$, the controller observes $x_t$ and chooses
$u_t=\pi_t(w_0,\ldots,w_{t-1})$. We assume each $\Xi_t$ has full column
rank, so past disturbances can be reconstructed from state transitions
and applied controls. Let $\Pi_{\mathrm c}$ be the causal policies whose
controls are square-integrable under every joint law in $\mathcal P_{\mathrm{drift}}$.
We restrict to the affine policy class
\begin{equation}\label{eq:drift_affine_class}
 \Pi_{\mathrm{aff}}^w=\left\{\pi\in\Pi_{\mathrm c}:
       u_t=g_t+\sum_{s<t}F_{ts}w_s\right\}.
\end{equation}

\section{SDP Reformulation}\label{sec:drift_sdp}

\subsection{Known-Law Controller}\label{subsec:drift_comparator}

We first compute the value and optimal policy of the comparator
$J^\star(\mathbb P)$. Fix a joint law $\mathbb P\in\mathcal P_{\mathrm{drift}}$
with moment path $(\boldsymbol\mu,\boldsymbol\Sigma)$, and let
$V_t:\R^n\to\R$ denote the minimum expected remaining cost as a
function of the state at time $t$.
Bellman's recursion gives
\begin{equation}\label{eq:drift_value_function}
 V_t(x)=x^\top S_tx+2x^\top P_t\boldsymbol\mu+c_t,
\end{equation}
where the coefficients $S_t\in\R^{n\times n}$,
$P_t\in\R^{n\times dT}$, and $c_t\in\R$ are computed backward from
$S_T=Q_T$, $P_T=0$, and $c_T=0$.
Let $E_t\in\R^{d\times dT}$ select the $t$th block of
$\boldsymbol\mu$. For $t=T-1,\ldots,0$, first compute
$M_t\in\R^{m\times m}$, $C_t,K_t\in\R^{m\times n}$, and $S_t$ by
\begin{align}
 M_t&=R_t+B_t^\top S_{t+1}B_t,\label{eq:drift_M}\\
 C_t&=B_t^\top S_{t+1}A_t+D_t^\top,\label{eq:drift_C}\\
 K_t&=-M_t^{-1}C_t,\label{eq:drift_K}\\
 S_t&=Q_t+A_t^\top S_{t+1}A_t-C_t^\top M_t^{-1}C_t.
 \label{eq:drift_S}
\end{align}
We assume $M_t\succ0$ for every stage. The matrices in
\eqref{eq:drift_M}--\eqref{eq:drift_S} are independent of the moment path.
Next, compute $P_t$ and the feedforward coefficient
$H_t\in\R^{m\times dT}$ from
\begin{align}
 P_t&=(A_t+B_tK_t)^\top
          (S_{t+1}\Xi_tE_t+P_{t+1}),\label{eq:drift_P}\\
 H_t&=-M_t^{-1}B_t^\top
          (S_{t+1}\Xi_tE_t+P_{t+1}).\label{eq:drift_H}
\end{align}
The scalar coefficient is updated by
\begin{equation}\label{eq:drift_c}
 c_t=c_{t+1}+\mu_t^\top\Xi_t^\top S_{t+1}\Xi_t\mu_t
 +2(P_{t+1}\boldsymbol\mu)^\top\Xi_t\mu_t
 +\tr(\Xi_t^\top S_{t+1}\Xi_t\Sigma_t)
 -(H_t\boldsymbol\mu)^\top M_t(H_t\boldsymbol\mu).
\end{equation}
The optimal causal control is
\begin{equation}\label{eq:drift_oracle}
 u_t^\star=K_tx_t+H_t\boldsymbol\mu.
\end{equation}
The comparator value is $J^\star(\mathbb P)=V_0(x_0)$.
The strict positivity of $M_t$ gives a unique minimizing action and a
finite value for every admissible moment path.

Set $\cH=\col(H_0,\ldots,H_{T-1})\in\R^{mT\times dT}$ and
$\mathcal M=\diag(M_0,\ldots,M_{T-1})\in\R^{mT\times mT}$.
The matrix $\cH$ is block upper triangular, with blocks
$\cH_{ts}\in\R^{m\times d}$.
The \emph{nominal controller} is the optimal controller when all stage
laws equal the reference law $\hat{\mathbb P}$, and is given by
\begin{equation}\label{eq:drift_nominal}
u_t^{\mathrm{nom}}=K_tx_t+\bar H_t\hat\mu,
\end{equation}
where $\bar H_t:=\sum_{s=t}^{T-1}\cH_{ts}\in\R^{m\times d}$ is the \emph{common-mean coefficient}.

\subsection{Derivation of Regret}\label{subsec:drift_regret}

To express regret compactly, we use the following equivalent
parameterization of the affine policy class \eqref{eq:drift_affine_class}:
\begin{equation}\label{eq:drift_residual}
 u_t=K_tx_t+\bar H_t\hat\mu
          +\tilde g_t+\sum_{s<t}\tilde F_{ts}(w_s-\hat\mu).
\end{equation}
Collect the coefficients as
$\tilde g=\col(\tilde g_0,\ldots,\tilde g_{T-1})\in\R^{mT}$ and
$\tilde F=(\tilde F_{ts})_{t,s=0}^{T-1}\in\R^{mT\times dT}$,
with $\tilde g_t\in\R^m$, $\tilde F_{ts}\in\R^{m\times d}$, and
$\tilde F_{ts}=0$ for $s\ge t$.
For this parameterization, write
$\mathfrak R(\tilde g,\tilde F)=\mathfrak R(\pi)$ and
$\mathfrak R(\tilde F):=\mathfrak R(0,\tilde F)$.
Define the \emph{mean deviations} $z_t=\mu_t-\hat\mu$ and
$z=\col(z_0,\ldots,z_{T-1})\in\R^{dT}$.
We next separate the mean and covariance contributions to the budgets.
Define $D_E,C_{\mathrm{sum}}\in\R^{dT\times dT}$ by
\[
 D_Ez=\col(z_0,z_1-z_0,\ldots,z_{T-1}-z_{T-2}),\qquad
 (C_{\mathrm{sum}}z)_t=\sum_{s=0}^t z_s.
\]
Define the positive definite matrices
\begin{equation}\label{eq:drift_budget_matrices}
 Q_E=D_E^\top\diag(a_0I_d,\ldots,a_{T-1}I_d)D_E,\quad
 Q_L=\diag(b_0I_d,\ldots,b_{T-1}I_d)
                                      +\eta C_{\mathrm{sum}}^\top C_{\mathrm{sum}}.
\end{equation}
Abbreviate the covariance contributions to the budgets by
\begin{align}
 e(\boldsymbol\Sigma)&=\sum_{t=0}^{T-1}a_t\mathsf B^2(\Sigma_t,\Sigma_{t-1}),
 \label{eq:drift_covariance_budgets}\\
 l(\boldsymbol\Sigma)&=\sum_{t=0}^{T-1}b_t\mathsf B^2(\Sigma_t,\hat\Sigma)
 +\eta\sum_{t=0}^{T-1}\mathsf B^2\!\left(
       \sum_{s=0}^t\Sigma_s,(t+1)\hat\Sigma\right).\notag
\end{align}

\begin{theorem}[Regret Representation]\label{thm:drift_regret}
Assume $M_t\succ0$ for all $t$. For a policy $\pi$ parameterized by
\eqref{eq:drift_residual} and a joint law $\mathbb P\in\mathcal P_{\mathrm{drift}}$
with moment path $(\boldsymbol\mu,\boldsymbol\Sigma)$, the regret depends
on $\mathbb P$ only through $z$ and $\boldsymbol\Sigma$ and is given by
\begin{equation}\label{eq:drift_full_regret}
 \mathcal R(\tilde g,\tilde F;z,\boldsymbol\Sigma)
   =\|\tilde g+(\tilde F-\cH)z\|_{\mathcal M}^2
   +\sum_{t=0}^{T-1}\tr(\mathcal A_t(\tilde F)\Sigma_t),
\end{equation}
where
\begin{equation}\label{eq:drift_coefficients}
 \mathcal A_t(\tilde F)=\tilde F_{:t}^\top\mathcal M\tilde F_{:t},
\end{equation}
with $\tilde F_{:t}\in\R^{mT\times d}$ denoting the $t$th block column.
$\mathfrak R(\tilde g,\tilde F)$ equals the attained optimal value of
\begin{equation}\label{eq:drift_full_worst_case_regret}
 \begin{aligned}
 \max_{z,\,\boldsymbol\Sigma}\quad
 &\|\tilde g+(\tilde F-\cH)z\|_{\mathcal M}^2
 +\sum_{t=0}^{T-1}\tr(\mathcal A_t(\tilde F)\Sigma_t)\\[0.3em]
 \text{s.t.}\quad
 &\Sigma_t\succeq0\quad(0\le t<T),\\
 &z^\top Q_Ez+e(\boldsymbol\Sigma)\le\delta_{\mathcal E}^2,\quad
 z^\top Q_Lz+l(\boldsymbol\Sigma)\le\delta_{\mathcal L}^2.
 \end{aligned}
\end{equation}
\end{theorem}

\begin{drproof}
Fix $\mathbb P$ and let $(x_t^\pi,u_t^\pi)$ denote the trajectory
of the policy $\pi$ whose regret is being evaluated.
With $\ell_t$ denoting the stage cost in \eqref{eq:drift_cost},
define the comparator's action-value function
$\mathsf Q_t^\star:\R^n\times\R^m\to\R$ by
\begin{equation}\label{eq:drift_action_value}
 \mathsf Q_t^\star(x,u)=\ell_t(x,u)
 +\E_{w_t\sim\mathbb P_t}\!\left[
 V_{t+1}(A_tx+B_tu+\Xi_tw_t)\right].
\end{equation}
The finite-horizon performance-difference identity
\cite{kakade2002approximately} gives
\begin{equation}\label{eq:drift_advantage_identity}
 \mathcal R(\pi;\mathbb P)=\E_{\mathbb P}^{\pi}\!\left[
 \sum_{t=0}^{T-1}\bigl(
 \mathsf Q_t^\star(x_t,u_t)-V_t(x_t)\bigr)\right]
\end{equation}
where the state $x_t$ is driven by the controls $u_t$ of the policy $\pi$.
Using the LQR recursion in Section~\ref{subsec:drift_comparator} and
completing the square yields
\begin{equation}\label{eq:drift_bellman_square}
 \mathsf Q_t^\star(x,u)-V_t(x)
       =\|u-K_tx-H_t\boldsymbol\mu\|_{M_t}^2.
\end{equation}
Set $\varepsilon_t=w_t-\mu_t$ and
$\varepsilon=\col(\varepsilon_0,\ldots,\varepsilon_{T-1})$.
Substituting \eqref{eq:drift_residual} and
\eqref{eq:drift_bellman_square} into \eqref{eq:drift_advantage_identity}
gives
\begin{equation}\label{eq:drift_expected_regret}
 \mathcal R(\pi;\mathbb P)=\E_{\mathbb P}\!\left[
 \|\tilde g+(\tilde F-\cH)z+\tilde F\varepsilon\|_{\mathcal M}^2
 \right].
\end{equation}
The centered disturbances $\varepsilon_t$ have zero mean. Independence implies
$\operatorname{Cov}(\varepsilon)=\diag(\Sigma_0,\ldots,\Sigma_{T-1})$,
so expanding \eqref{eq:drift_expected_regret} yields
\eqref{eq:drift_full_regret}.
Separating the mean and covariance terms in
\eqref{eq:drift_energy}--\eqref{eq:drift_reach} gives the constraints in
\eqref{eq:drift_full_worst_case_regret}. The supremum in
\eqref{eq:drift_worst_case_regret} is over laws, but regret depends only
on their moments. Independent, possibly degenerate Gaussian stage laws
realize any maximizing moment path. The feasible moment set is nonempty
since it contains the constant nominal path. It is closed by continuity
of the Gelbrich distance and bounded by the reach term $\mathcal L_1$,
hence compact. The continuous objective therefore attains its maximum.
\end{drproof}

Equations \eqref{eq:drift_advantage_identity} and
\eqref{eq:drift_bellman_square} allow us to interpret worst-case regret
minimization as a statistical prediction problem.
Indeed, for fixed $x_0$, any causal policy can be written as
$u_t=K_tx_t+\psi_t(w_0,\ldots,w_{t-1})$, where
$\psi_t:\R^{dt}\to\R^m$ may be nonlinear.
The regret is then
\begin{equation}\label{eq:drift_prediction_regret}
 \mathcal R(\pi;\mathbb P)
 =\sum_{t=0}^{T-1}\E_{\mathbb P}\!\left[
 \|\psi_t(w_0,\ldots,w_{t-1})-H_t\boldsymbol\mu\|_{M_t}^2\right].
\end{equation}
Regret minimization over $\Pi_{\mathrm c}$ therefore amounts to choosing
admissible functions $\psi_t$ to minimize
the worst-case weighted prediction error over
$\mathbb P\in\mathcal P_{\mathrm{drift}}$.
The prediction problem is shaped by the geometry of the control
problem: only $H_t\boldsymbol\mu$ needs to be predicted, rather than
the full mean path $\boldsymbol\mu$, and errors are measured in the
$M_t$-weighted norm.
Here we specialize to affine prediction rules.

\subsection{Controller SDP}\label{subsec:drift_synthesis}

We compute the
optimal worst-case regret and an optimal affine policy by solving an
SDP whose coefficients are obtained from a single LQR solve as described in Section~\ref{subsec:drift_comparator}.
\begin{theorem}[Affine DRRO Synthesis]\label{thm:drift_synthesis}
The following SDP, with $a_T=0$ and $U_T^{22}=0$, attains the optimal value of \eqref{eq:drift_drro}:
\begin{equation}\label{eq:drift_synthesis_sdp}
 \begin{aligned}
 \min_{\tilde F,\,\gamma_E,\,\gamma_L,\,\{U_t,V_t,\Lambda_t\}_{t=0}^{T-1}}\quad
 &\gamma_E\delta_{\mathcal E}^2+\gamma_L\delta_{\mathcal L}^2
       -\Bigl(\gamma_Ea_0+\gamma_L\sum_tb_t\Bigr)\tr\hat\Sigma\\
 &+\tr\!\left[\left(U_0^{22}+\sum_tV_t^{22}\right)\hat\Sigma\right]
 +\sum_t(t+1)\tr\!\left[
          (\Lambda_t^{22}-\gamma_L\eta I_d)\hat\Sigma\right]\\[0.5em]
 \text{s.t.}\quad
 &\gamma_E,\gamma_L\ge0,\quad\tilde F_{ts}=0,\ s\ge t,\\[0.3em]
 &U_t=\begin{pmatrix}
 U_t^{11}&-\gamma_Ea_tI_d\\
 -\gamma_Ea_tI_d&U_t^{22}
 \end{pmatrix}\succeq0,\quad
 V_t=\begin{pmatrix}
 V_t^{11}&-\gamma_Lb_tI_d\\
 -\gamma_Lb_tI_d&V_t^{22}
 \end{pmatrix}\succeq0,\\[0.3em]
 &\Lambda_t=\begin{pmatrix}
 \Lambda_t^{11}&-\gamma_L\eta I_d\\
 -\gamma_L\eta I_d&\Lambda_t^{22}
 \end{pmatrix}\succeq0,\quad 0\le t<T,\\[0.3em]
 &\begin{pmatrix}
 \gamma_EQ_E+\gamma_LQ_L&(\tilde F-\cH)^\top\\
 \tilde F-\cH&\mathcal M^{-1}
 \end{pmatrix}\succeq0,\\[0.3em]
 &{\renewcommand{\arraystretch}{1.2}\left(\begin{array}{c|c}
 \begin{gathered}
 \bigl[\gamma_E(a_t+a_{t+1})+\gamma_L(b_t+\eta(T-t))\bigr]I_d\\
 -U_t^{11}-U_{t+1}^{22}-V_t^{11}-\displaystyle\sum_{r=t}^{T-1}\Lambda_r^{11}
 \end{gathered}&\tilde F_{:t}^\top\\[0.3em]\hline
 \tilde F_{:t}&\mathcal M^{-1}
 \end{array}\right)}\succeq0,\quad 0\le t<T.
 \end{aligned}
\end{equation}
For any optimal $\tilde F$ from \eqref{eq:drift_synthesis_sdp}, a DRRO-optimal
affine controller is
\begin{equation}\label{eq:drift_optimal_policy}
 u_t^{\mathrm{DRRO}}=K_tx_t+\bar H_t\hat\mu
                  +\sum_{s<t}\tilde F_{ts}(w_s-\hat\mu).
\end{equation}
\end{theorem}

The correction to the nominal controller in \eqref{eq:drift_optimal_policy}
can assign different matrix weights $\tilde F_{ts}$ to past disturbances.
In the common-stage-law model \cite{fiechtner2026common}, the correction
instead weights all past disturbances equally and depends on them only
through the sample mean $\frac{1}{t}\sum_{s=0}^{t-1}w_s$ for $t\ge1$.
With a common law, all observations inform the same constant mean.
Here the controller must track a changing mean, so an observation's
relevance can depend on when it was obtained.

\begin{drproof}
Averaging \eqref{eq:drift_full_regret} at the feasible mean deviations
$z$ and $-z$, with the same covariance path, gives
$\mathfrak R(\tilde g,\tilde F)
\ge\tilde g^\top\mathcal M\tilde g+\mathfrak R(0,\tilde F)$ and so every optimal policy has $\tilde g=0$.
Thus, fix $\tilde g=0$ and some $\tilde F$ and consider the partial maximization over $z$ in
\eqref{eq:drift_full_worst_case_regret} for a fixed feasible covariance
path. With $\tilde g=0$, its objective becomes
$z^\top\mathcal B(\tilde F)z
+\sum_t\langle\mathcal A_t(\tilde F),\Sigma_t\rangle$, where
\begin{equation}\label{eq:drift_mean_coefficient}
 \mathcal B(\tilde F)=(\tilde F-\cH)^\top\mathcal M(\tilde F-\cH).
\end{equation}
We abbreviate $\mathcal B=\mathcal B(\tilde F)$ and
$\mathcal A_t=\mathcal A_t(\tilde F)$ and write $r_E=\delta_{\mathcal E}^2-e(\boldsymbol\Sigma)$
and $r_L=\delta_{\mathcal L}^2-l(\boldsymbol\Sigma)$ for the nonnegative
remaining budgets. We obtain the following exact convex relaxation for the partial maximization over $z$:
\begin{equation}\label{eq:drift_partial_mean_relaxation}
 \max_z\left\{z^\top\mathcal Bz:
                  z^\top Q_Ez\le r_E,\ z^\top Q_Lz\le r_L\right\}
 =\max_{Z\succeq0}
   \left\{\langle\mathcal B,Z\rangle:
                  \langle Q_E,Z\rangle\le r_E,\ \langle Q_L,Z\rangle\le r_L\right\}
\end{equation}
This can be justified as follows: If $r_E,r_L>0$, the SDP feasible set is compact since $Q_L\succ0$.
The rank bound \cite[Thm.~2.2]{pataki1998rank} therefore gives an
optimizer of rank at most one, proving exactness and attainment.
If either remaining budget is zero, $Q_E,Q_L\succ0$ force $z=0$
and $Z=0$, so both maxima are zero.
Consequently, $\mathfrak R(\tilde F)$ equals the optimal value of
\begin{equation}\label{eq:drift_mean_duality}
 \begin{aligned}
 \max_{Z,\,\boldsymbol\Sigma}\quad
 &\langle\mathcal B,Z\rangle+\sum_t\langle\mathcal A_t,\Sigma_t\rangle\\[0.3em]
 \text{s.t.}\quad
 &Z\succeq0,\quad\Sigma_t\succeq0\quad(0\le t<T),\\
 &\langle Q_E,Z\rangle+e(\boldsymbol\Sigma)\le\delta_{\mathcal E}^2,\quad
 \langle Q_L,Z\rangle+l(\boldsymbol\Sigma)\le\delta_{\mathcal L}^2.
 \end{aligned}
\end{equation}
The maximum is attained because the relaxation is exact and a maximizing
moment path from Theorem~\ref{thm:drift_regret} lifts to a feasible point
with $Z=zz^\top$ and the same objective value.

Strong duality holds for \eqref{eq:drift_mean_duality} and its Lagrange dual attains its optimum \cite[Sec.~5.2.3]{boyd2004convex} because
$\Sigma_t=\hat\Sigma$, $Z=\epsilon I_{dT}$ is a Slater point for
sufficiently small $\epsilon>0$, the objective is linear, and the budget
constraints are affine in $Z$ and convex in $\boldsymbol\Sigma$ by joint
convexity of the squared Bures distance
\cite[Sec.~3.1]{nguyen_mean_covariance_2025}.
For multipliers $\gamma_E,\gamma_L\ge0$, the Lagrangian is
\begin{equation}\label{eq:drift_lagrangian}
 \mathscr L(\gamma_E,\gamma_L,Z,\boldsymbol\Sigma)
  =\langle\mathcal B,Z\rangle+\sum_t\langle\mathcal A_t,\Sigma_t\rangle
  +\gamma_E(\delta_{\mathcal E}^2-\langle Q_E,Z\rangle-e(\boldsymbol\Sigma))
  +\gamma_L(\delta_{\mathcal L}^2-\langle Q_L,Z\rangle-l(\boldsymbol\Sigma)).
\end{equation}
Setting $\Gamma=\gamma_EQ_E+\gamma_LQ_L$ and maximizing the $Z$-dependent
part of the Lagrangian over $Z\succeq0$ gives
\begin{equation}\label{eq:drift_mean_elimination}
 \sup_{Z\succeq0}\langle\mathcal B-\Gamma,Z\rangle
 =\sup_z z^\top(\mathcal B-\Gamma)z
 =\begin{cases}
 0,&\Gamma\succeq\mathcal B,\\
 +\infty,&\text{otherwise}.
 \end{cases}
\end{equation}
Thus dualizing the two budget inequalities in
\eqref{eq:drift_full_worst_case_regret} with $\tilde g=0$ gives the same
dual function before and after lifting $zz^\top$ to $Z$.
Eliminating $Z$ gives the attained minimum
\begin{equation}\label{eq:drift_budget_dual}
 \begin{aligned}
 \min_{\gamma_E,\,\gamma_L}\quad
 &\gamma_E\delta_{\mathcal E}^2+\gamma_L\delta_{\mathcal L}^2
 +\sup_{\Sigma_t\succeq0}\Bigl\{
       \sum_t\langle\mathcal A_t,\Sigma_t\rangle-\gamma_Ee(\boldsymbol\Sigma)
              -\gamma_Ll(\boldsymbol\Sigma)\Bigr\}\\[0.3em]
 \text{s.t.}\quad
 &\gamma_E,\gamma_L\ge0,\quad
   \gamma_EQ_E+\gamma_LQ_L\succeq\mathcal B.
 \end{aligned}
\end{equation}

For the covariance supremum, we use the attained representation
\cite[Lem.~1]{nguyen_mean_covariance_2025} of the Bures distance:
\begin{equation}\label{eq:drift_covariance_lift}
 \begin{aligned}
 \mathsf B^2(S,U)=\min_Y\quad&\tr(S+U-2Y)\\[0.3em]
 \text{s.t.}\quad&\begin{pmatrix}S&Y\\Y^\top&U\end{pmatrix}\succeq0.
 \end{aligned}
\end{equation}
Applying \eqref{eq:drift_covariance_lift} to each Bures term introduces
$X_t,Y_t,\Theta_t\in\R^{d\times d}$ and the blocks
\begin{equation}\label{eq:drift_primal_blocks}
 \mathcal C_t^E=\begin{pmatrix}\Sigma_t&X_t\\X_t^\top&\Sigma_{t-1}\end{pmatrix},
 \qquad
 \mathcal C_t^L=\begin{pmatrix}\Sigma_t&Y_t\\Y_t^\top&\hat\Sigma\end{pmatrix},
 \qquad
 \mathcal C_t^P=\begin{pmatrix}
 \sum_{s=0}^{t}\Sigma_s&\Theta_t\\
 \Theta_t^\top&(t+1)\hat\Sigma
 \end{pmatrix}.
\end{equation}
Defining
\begin{equation}\label{eq:drift_cov_expenditure}
 \begin{aligned}
 \Phi_E&=\sum_ta_t\tr(\Sigma_t+\Sigma_{t-1}-2X_t),\\
 \Phi_L&=\sum_tb_t\tr(\Sigma_t+\hat\Sigma-2Y_t)
 +\eta\sum_t\tr\!\left(
       \sum_{s=0}^{t}\Sigma_s+(t+1)\hat\Sigma
                         -2\Theta_t\right),
 \end{aligned}
\end{equation}
the covariance supremum in \eqref{eq:drift_budget_dual} can be written,
using \eqref{eq:drift_covariance_lift}--\eqref{eq:drift_primal_blocks}
and $\gamma_E,\gamma_L\ge0$, as
\begin{equation}\label{eq:drift_covariance_sdp}
 \sup_{\boldsymbol\Sigma,\,\{X_t,Y_t,\Theta_t\}}\quad
 \sum_t\langle\mathcal A_t,\Sigma_t\rangle
                  -\gamma_E\Phi_E-\gamma_L\Phi_L
 \qquad
 \text{s.t.}\quad
 \mathcal C_t^E,\mathcal C_t^L,\mathcal C_t^P\succeq0,\quad 0\le t<T.
\end{equation}

The point $\Sigma_t=\hat\Sigma$, $X_t=Y_t=\Theta_t=0$ is strictly
feasible in \eqref{eq:drift_covariance_sdp}. At minimizing multipliers in
\eqref{eq:drift_budget_dual}, the value of \eqref{eq:drift_covariance_sdp}
is finite because \eqref{eq:drift_budget_dual} attains the finite value
$\mathfrak R(\tilde F)$.
For the blocks $U_t,V_t,\Lambda_t$ in \eqref{eq:drift_synthesis_sdp},
define the affine functions $\Omega_t\in\mathbb S^d$ and $\Psi\in\R$ by
\begin{subequations}\label{eq:drift_dual_cov_bound}
\begin{equation}\label{eq:drift_Omega}
 \Omega_t=\bigl[\gamma_E(a_t+a_{t+1})
                    +\gamma_L(b_t+\eta(T-t))\bigr]I_d
             -U_t^{11}-U_{t+1}^{22}-V_t^{11}
                    -\sum_{r=t}^{T-1}\Lambda_r^{11},
\end{equation}
\begin{equation}\label{eq:drift_Psi}
 \begin{aligned}
 \Psi={}&\gamma_E\delta_{\mathcal E}^2+\gamma_L\delta_{\mathcal L}^2
       -\Bigl(\gamma_Ea_0+\gamma_L\sum_tb_t\Bigr)\tr\hat\Sigma\\
 &+\tr\!\left[\left(U_0^{22}+\sum_tV_t^{22}\right)\hat\Sigma\right]
 +\sum_t(t+1)\tr\!\left[
          (\Lambda_t^{22}-\gamma_L\eta I_d)\hat\Sigma\right].
 \end{aligned}
\end{equation}
\end{subequations}
By conic duality \cite[Thm.~2.4.1(3.a)]{ben_tal_nemirovski2001},
substituting the dual of \eqref{eq:drift_covariance_sdp} into \eqref{eq:drift_budget_dual}
shows that $\mathfrak R(\tilde F)$ equals the attained minimum
\begin{equation}\label{eq:drift_fixed_dual}
 \begin{aligned}
 \min_{\gamma_E,\,\gamma_L,\,\{U_t,V_t,\Lambda_t\}}\quad&\Psi\\[0.3em]
 \text{s.t.}\quad&\gamma_E,\gamma_L\ge0,\quad
   U_t,V_t,\Lambda_t\succeq0,\quad 0\le t<T,\\
 &\gamma_EQ_E+\gamma_LQ_L\succeq\mathcal B,\quad
   \Omega_t\succeq\mathcal A_t,\quad 0\le t<T.
 \end{aligned}
\end{equation}
Adding the strictly causal matrix $\tilde F$ to the minimization
variables and taking Schur complements with
$\mathcal M^{-1}\succ0$ gives \eqref{eq:drift_synthesis_sdp}.

Finally, the nominal constant moment path is feasible, so
\begin{equation}\label{eq:drift_coercivity}
 \mathfrak R(\tilde F)
 \ge\sum_t\tr(\tilde F_{:t}^\top\mathcal M\tilde F_{:t}\hat\Sigma)
 \ge\lambda_{\min}(\mathcal M)\lambda_{\min}(\hat\Sigma)
                                      \|\tilde F\|_F^2.
\end{equation}
Thus $\mathfrak R$ is coercive on the closed space of strictly causal
matrices. Since it is also lower semicontinuous as the supremum of
continuous functions of $\tilde F$, it attains its minimum there.
Together with fixed-policy
dual attainment, this yields an optimizer of \eqref{eq:drift_synthesis_sdp}. Substituting the optimal $\tilde F$ and
$\tilde g=0$ into \eqref{eq:drift_residual} yields the controller
\eqref{eq:drift_optimal_policy}.
\end{drproof}

\begin{remark}[DRO Comparison]\label{rem:drift_dro}
In DRRO, every optimal policy has $\tilde g=0$, making the
mean-dependent part of the objective in
\eqref{eq:drift_full_worst_case_regret} homogeneous quadratic and enabling the exact
relaxation \eqref{eq:drift_partial_mean_relaxation}. 
On the other hand, for a DRO formulation of the problem
$\inf_{\pi\in\Pi_{\mathrm{aff}}^w}
\sup_{\mathbb P\in\mathcal P_{\mathrm{drift}}}J(\pi;\mathbb P)$,
the inner worst-case objective is also quadratic in $z$ and linear in
$\boldsymbol\Sigma$, but generally contains a linear term in $z$
whose coefficient, unlike in DRRO, need not vanish at optimality.
With two budgets, the resulting nonhomogeneous mean maximization can
have a duality gap. Thus, despite its additional comparator optimization,
DRRO admits an exact SDP here, whereas the analogous DRO relaxation
need not be exact.
With only the kinetic-energy constraint, however, the mean maximization
has an exact SDP relaxation even with a linear term, yielding an exact
DRO SDP.
\end{remark}

\section{Worst-Case Moments and Laws}\label{sec:drift_worst_case}

We derive worst-case laws for policies of the form
\eqref{eq:drift_residual} with $\tilde g=0$, which we call
\emph{centered policies}.
For the same reasons as in Remark~\ref{rem:drift_dro}, the mean SDP
relaxation need not be exact for $\tilde g\ne0$ in the setting with both constraints.
On the other hand, with only the $\mathcal E$ or $\mathcal L$ budget,
the mean maximization admits an exact SDP reformulation even for
$\tilde g\ne0$; see \cite{fiechtner2026common} for a similar argument.
For fixed $\tilde F$, abbreviate $\mathcal A_t=\mathcal A_t(\tilde F)$
and $\mathcal B=\mathcal B(\tilde F)$ as defined in
\eqref{eq:drift_coefficients} and \eqref{eq:drift_mean_coefficient}.

\begin{theorem}[Worst-Case Measure Recovery]\label{thm:drift_recovery}
Assume $M_t\succ0$ for all $t$ and fix a centered policy $\tilde F$.
If $\mathcal A_t=0$ for all $t$ and $\mathcal B=0$, every law in
$\mathcal P_{\mathrm{drift}}$ is worst-case. Otherwise, let
$(\gamma_E^\star,\gamma_L^\star,\{U_t^\star,V_t^\star,\Lambda_t^\star\})$ be any optimizer of
\eqref{eq:drift_fixed_dual}. At least one of $\gamma_E^\star,\gamma_L^\star$ is
positive. The worst-case covariance path is unique and satisfies
$\Sigma_t^\star\succ0$ for every $0\le t<T$.
If $\gamma_L^\star>0$, it is given by
\begin{equation}\label{eq:drift_covariance_recovery_L}
 \Sigma_t^\star=(\gamma_L^\star b_t)^2
 (V_t^{11,\star})^{-1}\hat\Sigma(V_t^{11,\star})^{-1}.
\end{equation}
If $\gamma_L^\star=0$, then $\gamma_E^\star>0$ and it is given recursively by
\begin{equation}\label{eq:drift_covariance_recovery_E}
 \Sigma_t^\star=(\gamma_E^\star a_t)^2
 (U_t^{11,\star})^{-1}\Sigma_{t-1}^\star(U_t^{11,\star})^{-1},
\end{equation}
with $\Sigma_{-1}^\star=\hat\Sigma$ and all inverses well defined in both cases.
Define
\begin{equation}\label{eq:drift_recovery_budgets}
 D_\mu=\gamma_E^\star Q_E+\gamma_L^\star Q_L-\mathcal B,\qquad
 r_E=\delta_{\mathcal E}^2-e(\boldsymbol\Sigma^\star),\qquad
 r_L=\delta_{\mathcal L}^2-l(\boldsymbol\Sigma^\star).
\end{equation}
The worst-case mean deviations are precisely the vectors $z$ satisfying
\begin{equation}\label{eq:drift_mean_recovery}
 D_\mu z=0,\quad z^\top Q_Ez\le r_E,\quad z^\top Q_Lz\le r_L,\quad
 \gamma_E^\star(r_E-z^\top Q_Ez)=0,\quad
   \gamma_L^\star(r_L-z^\top Q_Lz)=0.
\end{equation}
If $\ker D_\mu=\{0\}$ or either remaining budget is zero, the only
solution is $z=0$. Otherwise, choose a basis matrix
$N\in\R^{dT\times k}$ of $\ker D_\mu$ and an orthogonal matrix $O$ such that
\begin{equation}\label{eq:drift_mean_basis}
 N^\top Q_EN=I_k,\qquad
 N^\top Q_LN=O\diag(\kappa_1,\ldots,\kappa_k)O^\top.
\end{equation}
If both multipliers are positive, all solutions are
\begin{equation}\label{eq:drift_mean_parametrization}
 z=\sqrt{r_E}\,NO\col(s_1\sqrt{p_1},\ldots,s_k\sqrt{p_k}),
\end{equation}
where $s_j\in\{-1,1\}$ and
\begin{equation}\label{eq:drift_mean_weights}
 p_j\ge0,\qquad \sum_jp_j=1,\qquad
 \sum_j\kappa_jp_j=\frac{r_L}{r_E}.
\end{equation}
If only $\gamma_E^\star>0$, replace the weighted-sum equality in
\eqref{eq:drift_mean_weights} by $\sum_j\kappa_jp_j\le r_L/r_E$.
If only $\gamma_L^\star>0$, interchange $E$ and $L$ in
\eqref{eq:drift_mean_basis}--\eqref{eq:drift_mean_weights} and use the
weighted-sum inequality.
For every such $z$, any independent stage laws with means
$\hat\mu+z_t$ and covariances $\Sigma_t^\star$ are worst-case.
In particular, one may take
$\mathbb P^\star=\bigotimes_{t=0}^{T-1}
\mathcal N(\hat\mu+z_t,\Sigma_t^\star)$.
The dual variables from an optimizer of \eqref{eq:drift_synthesis_sdp}
are optimal in \eqref{eq:drift_fixed_dual} for its returned policy.
\end{theorem}

\begin{drproof}
If both multipliers vanish, feasibility in \eqref{eq:drift_fixed_dual}
forces $\mathcal B=0$ and $\mathcal A_t=0$ for all $t$.
Thus at least one multiplier is positive in the nonconstant case.
Theorem~\ref{thm:drift_regret} and the proof of
Theorem~\ref{thm:drift_synthesis} give primal attainment and exact
duality with \eqref{eq:drift_fixed_dual}.
Take any worst-case moment path and form the blocks
$\mathcal C_t^E,\mathcal C_t^L$ in \eqref{eq:drift_primal_blocks}
using minimizing couplings from \eqref{eq:drift_covariance_lift}.
With $\gamma_E=\gamma_E^\star$ and $\gamma_L=\gamma_L^\star$,
complementary slackness between the covariance SDP
\eqref{eq:drift_covariance_sdp} and its dual used in
\eqref{eq:drift_fixed_dual} gives
$U_t^\star\mathcal C_t^E=0$ and $V_t^\star\mathcal C_t^L=0$.
If $\gamma_L^\star>0$, the nonsingular off-diagonal block of $V_t^\star\succeq0$
implies $V_t^{11,\star}\succ0$. The top-right and top-left blocks of
$V_t^\star\mathcal C_t^L=0$ give
$V_t^{11,\star}Y_t=\gamma_L^\star b_t\hat\Sigma$ and
$V_t^{11,\star}\Sigma_t=\gamma_L^\star b_tY_t^\top$,
which yield \eqref{eq:drift_covariance_recovery_L}.
If $\gamma_L^\star=0$, then $\gamma_E^\star>0$, so the same argument
with $U_t^\star$ gives
\eqref{eq:drift_covariance_recovery_E}.
Starting from $\hat\Sigma\succ0$, both formulas give positive definite
covariances. By strong duality, any two worst-case covariance paths
satisfy $U_t^\star\mathcal C_t^E=0$ and $V_t^\star\mathcal C_t^L=0$ with the same chosen dual
optimizer. The formula \eqref{eq:drift_covariance_recovery_L} or the
recursion \eqref{eq:drift_covariance_recovery_E} therefore forces the
paths to coincide.
At the chosen dual optimizer, strong duality and
\eqref{eq:drift_mean_elimination} give $z^\top D_\mu z=0$,
equivalently $D_\mu z=0$ since $D_\mu\succeq0$.
Together with feasibility and complementary slackness for the budget
constraints in \eqref{eq:drift_mean_duality}, this gives
\eqref{eq:drift_mean_recovery}.
Conversely, the recovered covariance path maximizes the covariance
supremum in \eqref{eq:drift_budget_dual}. Every $z$ satisfying
\eqref{eq:drift_mean_recovery} is feasible, maximizes the mean
Lagrangian in \eqref{eq:drift_mean_elimination}, and makes the budget
penalties vanish. Its regret therefore equals the optimal value of \eqref{eq:drift_fixed_dual}.

Writing $z=NOv$ gives
$z^\top Q_Ez=\|v\|^2$ and $z^\top Q_Lz=\sum_j\kappa_jv_j^2$.
When $\gamma_E^\star>0$, \eqref{eq:drift_mean_recovery} fixes
$\|v\|^2=r_E$, so write $v_j=\sqrt{r_E}s_j\sqrt{p_j}$.
The remaining budget gives \eqref{eq:drift_mean_weights}, with an
inequality when $\gamma_L^\star=0$. Interchanging the budgets covers
$\gamma_E^\star=0$. If either remaining budget is zero, $Q_E,Q_L\succ0$
force $z=0$.
The assertion about the laws follows because regret depends only on the moments.
Finally, an optimizer of \eqref{eq:drift_synthesis_sdp} must also
minimize \eqref{eq:drift_fixed_dual} at its returned policy.
\end{drproof}

\section{Experiments}\label{sec:drift_experiments}

We use the single-asset Almgren--Chriss portfolio liquidation model and
parameter values of \cite[Sec.~V]{fiechtner2026common}. With remaining inventory
$q_t$, sales $u_t$, and midprice $s_t$, we maximize the risk-adjusted
expected liquidation proceeds over $T=20$ trades,
\begin{equation}\label{eq:drift_execution_proceeds}
 J(\pi;\mathbb P):=
 \E_{\mathbb P}^{\pi}\!\left[\sum_{t=0}^{T-1}
 \left(u_t(s_t-Ru_t)-\lambda\hat\sigma^2q_t^2\right)\right],
\end{equation}
subject to complete liquidation, $q_T=0$, and the dynamics
\begin{equation}\label{eq:drift_execution_dynamics}
 q_{t+1}=q_t-u_t,\quad 0\le t\le T-1,
 \qquad
 s_{t+1}=s_t-Gu_t+w_t,\quad 0\le t\le T-2.
\end{equation}
Here $w_t$ is a price innovation with uncertain mean and variance.
Temporary impact reduces a sale's execution price to $s_t-Ru_t$,
whereas permanent impact lowers subsequent midprices by $Gu_t$.
The objective maximizes expected net trading proceeds less a quadratic
inventory penalty, with risk aversion $\lambda$ and fixed nominal
innovation variance $\hat\sigma^2$.
Purchases and short positions are permitted.
Inventories and trades are in millions of shares and prices in
dollars per share. We use $q_0=1$, $s_0=100$, $G=0.08$, $R=0.15$,
$\lambda=2$, $\hat\mu=0$, and $\hat\sigma=0.12$.
We eliminate the final trade using $u_{T-1}=q_{T-1}$, minimize $-J$,
and report $J$ in millions of dollars.
We solve the SDPs with Clarabel and compute expectations directly
from the moments.
For the $T-1$ relevant innovations, we set $a_t=1$, $b_t=1/(T-1)$,
$\eta=1/2470$, $\delta_{\mathcal E}=0.10$, and $\delta_{\mathcal L}=\sqrt{2}\,0.10$.
The common-law moment set is
\begin{equation}\label{eq:drift_common_moment_set}
 \cW_{\mathrm{com}}:=\left\{(\mu,\sigma^2):\sigma\ge0,\quad
 \mu^2+(\sigma-\hat\sigma)^2\le0.1^2\right\}.
\end{equation}
On constant moment paths, $\mathcal E=\mu^2+(\sigma-\hat\sigma)^2$
and $\mathcal L\le2\mathcal E$. Hence $\cW_{\mathrm{drift}}$ contains
every constant-path extension of $\cW_{\mathrm{com}}$, and its
restriction to constant moment paths is exactly
\eqref{eq:drift_common_moment_set}. The larger set imposes a stronger
robustness requirement: for any fixed policy, worst-case regret over
$\cW_{\mathrm{drift}}$ is at least that over $\cW_{\mathrm{com}}$.

We compare three policies: Nominal, DRRO over $\cW_{\mathrm{com}}$
(Common DRRO), and DRRO over $\cW_{\mathrm{drift}}$ (Drifting DRRO).
Table~\ref{tab:drift_execution} compares the policies on each of the
two evaluation sets separately.
On $\cW_{\mathrm{drift}}$, Drifting DRRO reduces worst-case
regret by \$1.13 million relative to Nominal and \$1.24 million
relative to Common DRRO.
This protection comes at a cost in the settings for which the
baselines are optimized: the nominal objective value is \$0.13 million lower than
Nominal's, and worst-case regret over $\cW_{\mathrm{com}}$ is
\$0.38 million higher than Common DRRO's.
Both performance losses are substantially smaller than the corresponding
reductions in worst-case regret over $\cW_{\mathrm{drift}}$.

\begin{table}[!htbp]
\centering
\caption{Liquidation performance in millions of dollars.}
\label{tab:drift_execution}
\setlength{\tabcolsep}{7pt}
\begin{tabular}{lrrr}
\hline
Policy & Nominal objective value & \multicolumn{2}{c}{Worst-case regret}\\
 & & $\cW_{\mathrm{com}}$ & $\cW_{\mathrm{drift}}$\\
\hline
Nominal & \DriftNominalReturn & \DriftNominalCommonRegret & \DriftNominalRegret\\
Common DRRO & \DriftCommonReturn & \DriftCommonCommonRegret & \DriftCommonRegret\\
Drifting DRRO & \DriftRobustReturn & \DriftRobustCommonRegret & \DriftRobustRegret\\
\hline
\end{tabular}
\end{table}

\begin{figure}[!htbp]
\centering
\includegraphics[width=0.72\linewidth]{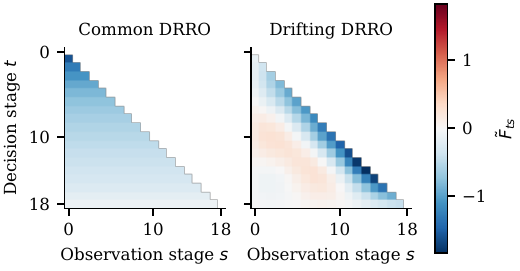}
\caption{Residual feedback coefficients $\tilde F_{ts}$.}
\label{fig:drift_feedback}
\end{figure}

Figure~\ref{fig:drift_feedback} shows the coefficients $\tilde F_{ts}$ from
\eqref{eq:drift_optimal_policy}. Common DRRO weights past innovations
equally at each trading time, whereas Drifting DRRO assigns the
strongest negative weights to the most recent observations.
Holding the current state fixed, positive recent innovations therefore
reduce sales and may induce purchases.
Positive weights on some older observations instead increase sales after
earlier price rises, consistent with the reach-and-persistence budget
limiting the continuation of sustained upward mean shifts.

\begin{figure}[!htbp]
\centering
\includegraphics[width=\textwidth]{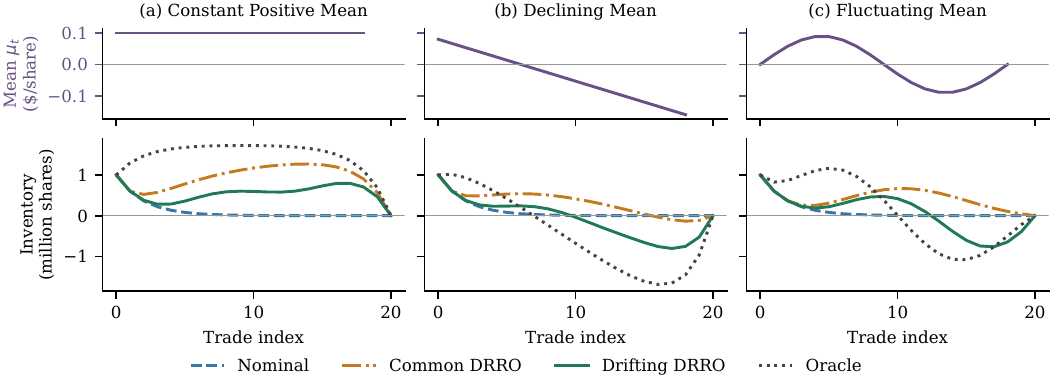}
\caption{Prescribed innovation means (upper) and expected inventories
(lower), with $\Sigma_t=\hat\sigma^2$ throughout. The mean paths are
(a) $\mu_t=0.10$, (b) $\mu_t=0.08-0.24t/(T-2)$, and
(c) $\mu_t=0.09\sin(2\pi t/(T-2))$.
All controllers face the same moment path within each column.}
\label{fig:drift_inventory}
\end{figure}

Figure~\ref{fig:drift_inventory} compares expected inventories under
three prescribed moment paths in $\cW_{\mathrm{drift}}$, with innovation
variance fixed at its reference value $\hat\sigma^2$.
The constant positive-mean case also belongs to $\cW_{\mathrm{com}}$.
Oracle uses the known-law optimal causal controller
\eqref{eq:drift_oracle}.
Common DRRO follows Oracle's expected inventory more closely overall
in the constant-mean case, whereas Drifting DRRO does so in the
declining- and fluctuating-mean cases.
The same ordering holds for regret.

\section{Conclusion}\label{sec:drift_conclusion}

Under the proposed drift budgets, affine ex-ante DRRO admits an exact
SDP that yields both an optimal controller and worst-case laws.
Liquidation experiments demonstrate substantially lower worst-case
regret over the drifting ambiguity set than common-law DRRO and nominal
control.
Partial observability is a natural extension.

\section*{Acknowledgments}
We acknowledge the use of GPT 6 Astra, OpenAI, \cite{chatgpt2026}, for editorial assistance, writing code for experiments,
and verifying mathematical derivations. We reviewed and edited all outputs and take full responsibility
for the scientific integrity of the final content.
J. Blanchet also gratefully acknowledges support from ONR under award N00014-24-1-2655, and the National Science Foundation (NSF) under grants 2312204 and 2403007.

\bibliographystyle{IEEEtran}
\bibliography{references}
\end{document}